\documentclass[runningheads]{llncs}
\usepackage{cite}
\usepackage{amsmath,amssymb,amsfonts}
\usepackage{algorithm,algorithmic}
\usepackage[dvipdfm]{graphicx}
\usepackage{url}
\usepackage{textcomp}
\usepackage{xcolor}
\def\BibTeX{{\rm B\kern-.05em{\sc i\kern-.025em b}\kern-.08em
    T\kern-.1667em\lower.7ex\hbox{E}\kern-.125emX}}
\begin{document}

\title{Acceleration of multiple precision matrix multiplication based on multi-component floating-point arithmetic using AVX2}

\author{Tomonori Kouya\inst{1}\orcidID{0000-0003-0178-5519}}
\authorrunning{T.Kouya}
%
\institute{Shizuoka Institute of Science and Technology,\\2200-2 Toyosawa, Fukuroi 437-8555, Japan \\
\email{kouya.tomonori@sist.ac.jp}}

\maketitle

\begin{abstract}
In this paper, we report the results obtained from the acceleration of multi-binary64-type multiple precision matrix multiplication with AVX2. We target double-double (DD), triple-double (TD), and quad-double (QD) precision arithmetic designed by certain types of error-free transformation (EFT) arithmetic. Furthermore, we implement SIMDized EFT functions, which simultaneously compute with four binary64 numbers on x86\_64 computing environment, and by using help of them, we also develop SIMDized DD, TD, and QD additions and multiplications. In addition, AVX2 load/store functions were adopted to efficiently speed up reading and storing matrix elements from/to memory. Owing to these combined techniques, our implemented multiple precision matrix multiplications have been accelerated more than three times compared with non-accelerated ones. Our accelerated matrix multiplication modifies the performance of parallelization with OpenMP.
\keywords{multiple precision floating-point arithmetic \and matrix multiplication \and AVX2 \and parallelization}
\end{abstract}

%
\section{Introduction}

Scale of scientific computation is becoming increasingly larger, and to obtain the user-required accuracy of approximation for ill-conditioned cases, it requires more digits of fraction in floating-point number than those of the IEEE754 binary standard, such as binary32 (23-bits fraction) and binary64 (53-bits). To obtain precise approximation, it is critical to adopt the multiple precision floating-point (MPF) arithmetic, which can extend its digits of fraction more than that of binary64.

The implementation of MPF arithmetic is currently categorized into two groups: multi-digit-type based on integer arithmetic and multi-component-type (or multi-term), which is constructed by several binary32 or binary64 numbers. The QD library by Bailey et al.\cite{qd} represents the latter, whereas MPFR\cite{mpfr} based on the MPN kernel of GNU MP\cite{gmp} represents the former. Both libraries are coded by C and C++ languages, and the MPN kernel is accelerated by assembler codes for various CPU architectures. 

MPLAPACK based on MPBLAS\cite{mplapack} is an efficient MPF linear computation library, which employs QD, MPFR, and bridge functions. Although it is implemented from Fortran codes of reference BLAS, more accelerated MPF BLAS libraries, such as ATLAS\cite{atlas} or OpenBLAS\cite{openblas}, are expected.

It is well-known that multi-component MPF arithmetic can be accelerated using SIMD instructions, such as AVX2. Lis\cite{lis} is one of the most optimal AVX2-based linear and sparse linear computation libraries, which has achieved over three times speedup using AVX2\cite{hasegawa_mupat}. 

In this paper, we report the performance of our implemented binary64-based multi-component MPF matrix multiplication using a technique similar to that of Lis. The targeted precisions are double-double (DD), triple-double (TD), and quadruple-double (QD) arithmetic. We adopted SIMDized error-free transformation functions on the x86\_64 CPU architecture to extend the performance of these three precision MPF linear computations. QD-based addition is also utilized to obtain better performance in the TD arithmetic. A new data structure for the MPF vector and matrix is introduced to accelerate the load and store their elements with AVX2. Owing to these efforts, we have realized higher performance for all types of MPF arithmetic than those of the Rgemm function of MPBLAS. Furthermore, we have modified the speed of parallelization with OpenMP for block matrix multiplication in the AMD EPYC environment.

We determined that our Strassen matrix multiplication could not be accelerated on a CPU with over eight cores and that block matrix multiplication is better performed on them.

%
\section{Acceleration of DD, TD, and QD matrix multiplication with AVX2}

Multi-component MPF arithmetic are constructed by several hard-wired binary32 or binary64 numbers and various functions of the so-called error-free transformation (EFT) \cite{dekker}. The DD, TD, and QD precision MPF numbers are constructed with two binary64, three binary64, and four binary64 numbers, respectively. The most significant component (one leading binary64 number) of DD, TD, and QD precision numbers possess the most significant digits of the entire MPF number. 

It has been reported that, using DD precision arithmetic, basic linear computation can be accelerated with SIMD instructions, including AVX2, as DD arithmetic does not include complicated branches in its addition, subtraction, multiplication, and division. In contrast, TD and QD arithmetic require renormalization processes, including some conditional branches; therefore, we are required to implement these renormalization as simple loops, where we consider each binary64 number in the \verb|_m256d| data type of AVX2, including three or four binary64 numbers. Because we completed all DD, TD, and QD arithmetic using the \verb|_m256d| data type as API interfaces, our MPF matrix multiplication library has been simply implemented and accelerated with AVX2.

%
\subsection{Error-free transformation with AVX2}

As we explained earlier, we only implemented the \verb|_m256d| data type to express and accelerate linear computation with multi-component MPF arithmetic. Therefore, C APIs such as \verb|_mm256_[add, sub, mul, div]_pd|, and \verb|_mm256_fmadd_pd| as fused multiply add (FMA) \cite{intel_avx} were embedded in required EFT functions, such as QuickTwoSum, TwoSum, and TwoProd. Subsequently, these EFT functions were implemented as AVX2QuickTwoSum, AVX2TwoSum, and AVX2TwoProd, respectively.
    
%
\subsection{SIMDized DD addition and multiplication}

As DD arithmetic is the simplest of the multi-component MPF arithmetic, its addition and multiplication needed for matrix multiplication, which directly adopt SIMDized EFT functions, can be easily implemented, as shown in Algorithm \ref{algo:AVX2DDadd} and \ref{algo:AVX2DDmul}. One DD precision number is constructed with two binary64 as $x[2] = (x[0], x[1])$, and each array element is expressed as one \verb|_m256d| variable. AVX2DDadd and AVX2DDmul simultaneously execute four calculations with $x[2]$ and $y[2]$, in which they have four DD numbers.

\begin{algorithm}
    \caption{ $r[2] :=\mathrm{AVX2DDadd}(x[2], y[2])$ }\label{algo:AVX2DDadd}
    \begin{algorithmic}
        \STATE $(s, e) := \mbox{AVX2TwoSum}(x[0], y[0])$
        \STATE $w := $\verb|_mm256_add_pd|($x[1]$, $y[1]$)
        \STATE $e := $\verb|_mm256_add_pd|($e$, $w$)
        \STATE $(r[0], r[1]) := \mbox{AVX2QuickTwoSum}(s, e)$
        \STATE \textbf{return} ($r[0]$, $r[1]$)
    \end{algorithmic}
\end{algorithm}

\begin{algorithm}
    \caption{ $r[2] :=\mathrm{AVX2DDmul}(x[2], y[2])$ }\label{algo:AVX2DDmul}
    \begin{algorithmic}
        \STATE $(p_1, p_2) := \mbox{AVX2TwoProd}(x[0], y[0])$
        \STATE $w_1 := $\verb|_mm256_mul_pd|($x[0], y[1]$)
        \STATE $w_2 := $\verb|_mm256_mul_pd|($x[1], y[0]$)
        \STATE $w_3 := $\verb|_mm256_add_pd|($w_1, w_2$)
        \STATE $p_2 := $\verb|_mm256_add_pd|($p_2, w_3$) 
        \STATE $(r[0], r[1]) :=\mbox{AVX2QuickTwoSum}(p_1, p_2)$
    \end{algorithmic}
\end{algorithm}

As previously introduced, Lis\cite{lis} and MuPaT\cite{hasegawa_mupat} realized a DD linear computation with AVX2 that is approximately three times faster. They utilize C macros to implement the SIMDized DD arithmetic. In contrast to Lis and MuPat, we adopted C static inline functions of SIMDized EFT for DD, TD, and QD linear computations. In the remainder of this section, we describe these algorithms with AVX2 and their performance.

%
\subsection{SIMDized TD addition and multiplication}

Fabiano et.al. proposed the optimized triple word arithmetic\cite{tripleword}. In this arithmetic set, the renormalization of triple word number adopts VecSum and VSEB($k$) (VecSum with Blanch). Our SIMDized TD arithmetic can be implemented with AVX2 instruction, AVX2VecSum and AVX2VSEB($n$), respectively. AVX2VecSum can be fully implemented with AVX2 functions, whereas AVX2VSEB($n$) is not completely SIMDized owing to the conditional branch embedded in it.

TD addition (TDadd) calculates $r[3] = (r[0], r[1], r[2])$ as the sum of $x[3] = (x[0], x[1], x[2])$ and $y[3] = (y[0], y[1], y[2])$. The standard approach to obtain it is by merge sorting all elements of $x[3]$ and $y[3]$, normalizing them by VecSum, and expressing $r[3]$ as the TD number by VSEB(3). The following AVX2TDadd (Algorithm \ref{algo:AVX2TDsum}) is implemented with AVX2. This algorithm is not sufficiently accelerated with AVX2 owing to the poor performance of the merge(-sort) function and incomplete AVXVSEB($n$).

\begin{algorithm}
    \caption{ $r[3] :=$ AVX2TDadd($x[3]$, $y[3]$)}\label{algo:AVX2TDsum}
    \begin{algorithmic}
        \STATE $(z_0, ..., z_5) := \mbox{AVX2Merge}(x[0], x[1], x[2], y[0], y[1], y[2])$
        \STATE $(e_0, ..., e_5) := \mbox{AVX2VecSum}(z_0, ..., z_5)$
        \STATE $(r[0], r[1], r[2]) := \mbox{AVX2VSEB(3)}(e_0, ..., e_5)$
        \STATE \textbf{return} ($r[0]$, $r[1]$, $r[2]$)
    \end{algorithmic}
\end{algorithm}

As demonstrated in the results obtained by benchmark tests later, AVX2TDadd cannot be accelerated completely with AVX2. Therefore, as an alternative, we applied QDadd (Algorithm \ref{algo:AVX2QDadd}) as TD addition (TDaddq) with a substitution such as $x[3] = y[3] = 0$. In our TD matrix multiplication, we adopted TDaddq and AVX2TDaddq as default TD additions.

For the TD multiplication, Fabiano et al. proposed a heavy and precise ``accurate'' version, and a light and sloppy ``fast'' version. We selected the latter as the default TD multiplication (TDmul) and implemented AVX2TDmul with AVX2.

\begin{algorithm}
    \caption{$r[3]$ $:=$ AVX2TDmul($x[3]$, $y[3]$)}\label{algo:AVX2TDprodSloppy}
    \begin{algorithmic}
        \STATE $(z_{00}^{\rm up}, z_{00}^{\rm lo}) := \mbox{AVX2TwoProd}(x[0], y[0])$
        \STATE $(z_{01}^{\rm up}, z_{01}^{\rm lo}) := \mbox{AVX2TwoProd}(x[0], y[1])$
        \STATE $(z_{10}^{\rm up}, z_{10}^{\rm lo}) := \mbox{AVX2TwoProd}(x[1], y[0])$
        \STATE $(b_0, b_1, b_2) := \mbox{AVX2VecSum}(z_{00}^{\rm lo}, z_{01}^{\rm up}, z_{10}^{\rm up})$
        \STATE $c := $\verb|_mm256_fmadd_pd(|$x[1], y[1], b_2$\verb|)|
        \STATE $z_{31} := $\verb|_mm256_fmadd_pd(|($x[0], y[2], z_{10}^{\rm lo}$\verb|)|
        \STATE $z_{32} := $\verb|_mm256_fmadd_pd(| $x[2], y[0], z_{01}^{\rm lo}$ \verb|)|
        \STATE $z_{3} :=$ \verb|_mm256_add_pd(|$z_{31}$, $z_{32}$\verb|)|
        \STATE $s_{3} :=$ \verb|_mm256_add_pd(|$c$, $z_{3}$\verb|)|
        \STATE $(e_0, e_1, e_2, e_3) := \mbox{AVX2VecSum}(z_{00}^{\rm up}, b_0, b_1, s_3)$
        \STATE $r[0] := e_0$
        \STATE $(r[1], r[2]) := \mbox{AVX2VSEB}(2)(e_1, e_2, e_3)$
        \STATE \textbf{return} ($r[0]$, $r[1]$, $r[2]$)
    \end{algorithmic}
\end{algorithm}

As a trial approach, TDmulq, which was converted to QDmul for the TD multiplication, was implemented, although it could not realize a higher speed than TDmul. Therefore, to implement the TD matrix multiplication, we adopted the TDaddq and TDmul sets, as well as those of AVX2TDaddq and AVX2TDmul.

%
\subsection{SIMDized QD addition and multiplication}

For QD arithmetic, we selected a lighter addition and multiplication (sloppy addition and multiplication), and then we implemented AVX2QDadd (Algorithm\ \ref{algo:AVX2QDadd}) and AVX2QDmul(Algorithm \ref{algo:AVX2QDmul}) with AVX2Renorm as the renormalization of QD numbers with AVX2.

\begin{algorithm}
    \caption{ $r[4] :=\mathrm{AVX2QDadd}(x[4], y[4])$ }\label{algo:AVX2QDadd}
    \begin{algorithmic}
        \STATE $s_0 := $\verb|_mm256_add_pd|($x[0], y[0]$)
        \STATE $s_1 := $\verb|_mm256_add_pd|$(x[1], y[1])$
        \STATE $s_2 := $\verb|_mm256_add_pd|$(x[2], y[2])$
        \STATE $s_3 := $\verb|_mm256_add_pd|$(x[3], y[3])$
        \STATE $v_0 := $\verb|_mm256_sub_pd|$(s_0, x[0])$
        \STATE $v_1 := $\verb|_mm256_sub_pd|$(s_1, x[1])$
        \STATE $v_2 := $\verb|_mm256_sub_pd|$(s_2, x[2])$
        \STATE $v_3 := $\verb|_mm256_sub_pd|$(s_3, x[3])$
        \STATE $u_0 := $\verb|_mm256_sub_pd|$(s_0, v_0)$
        \STATE $u_1 := $\verb|_mm256_sub_pd|$(s_1, v_1)$
        \STATE $u_2 := $\verb|_mm256_sub_pd|$(s_2, v_2)$
        \STATE $u_3 := $\verb|_mm256_sub_pd|$(s_3, v_3)$
        \STATE $w_0 := $\verb|_mm256_sub_pd|$(x[0], u_0)$
        \STATE $w_1 := $\verb|_mm256_sub_pd|$(x[1], u_1)$
        \STATE $w_2 := $\verb|_mm256_sub_pd|$(x[2], u_2)$
        \STATE $w_3 := $\verb|_mm256_sub_pd|$(x[3], u_3)$
        \STATE $u_0 := $\verb|_mm256_sub_pd|$(y[0], v_0)$
        \STATE $u_1 := $\verb|_mm256_sub_pd|$(y[1], v_1)$
        \STATE $u_2 := $\verb|_mm256_sub_pd|$(y[2], v_2)$
        \STATE $u_3 := $\verb|_mm256_sub_pd|$(y[3], v_3)$
        \STATE $t_0 := $\verb|_mm256_add_pd|$(w_0, u_0)$
        \STATE $t_1 := $\verb|_mm256_add_pd|$(w_1, u_1)$
        \STATE $t_2 := $\verb|_mm256_add_pd|$(w_2, u_2)$
        \STATE $(s_1, t_0) := \mathrm{AVX2TwoSum}(s_1, t_0)$
        \STATE $(s_2, t_0, t_1) := \mathrm{AVX2ThreeSum}(s_2, t_0, t_1)$
        \STATE $(s_3, t_0) := \mathrm{AVX2ThreeSum2}(s_3, t_0, t_2)$
        \STATE $t_0 := $\verb|_mm256_add_pd|(\verb|_mm256_add_pd|($t_0, t_1$), $t_3$)
        \STATE $(r[0], r[1], r[2], r[3]) := \mathrm{AVX2Renorm}(s_0, s_1, s_2, s_3, t_0)$
        \STATE \textbf{return} ($r[0], r[1], r[2], r[3]$)
    \end{algorithmic}
\end{algorithm}

\begin{algorithm}
    \caption{ $r[4] :=\mathrm{AVX2QDmul}(x[4], y[4])$ }\label{algo:AVX2QDmul}
    \begin{algorithmic}
        \STATE $s_0 := $\verb|_mm256_add_pd|($x[0]$, $y[0]$)
        \STATE $(p_0, q_0) := \mathrm{AVX2TwoProd}(x[0], y[0])$
        \STATE $(p_1, q_1) := \mathrm{AVX2TwoProd}(x[0], y[1])$
        \STATE $(p_2, q_2) := \mathrm{AVX2TwoProd}(x[1], y[0])$
        \STATE $(p_3, q_3) := \mathrm{AVX2TwoProd}(x[0], y[2])$
        \STATE $(p_4, q_4) := \mathrm{AVX2TwoProd}(x[1], y[1])$
        \STATE $(p_5, q_5) := \mathrm{AVX2TwoProd}(x[2], y[0])$
        \STATE $(p_1, p_2, q_0) := \mathrm{AVX2ThreeSum}(p_1, p_2, q_0)$
        \STATE $(p_2, q_1, q_2) := \mathrm{AVX2ThreeSum}(p_2, q_1, q_2)$
        \STATE $(p_3, p_4, p_5) := \mathrm{AVX2ThreeSum}(p_3, p_4, p_5)$
        \STATE $(s_0, t_0) := \mathrm{AVX2TwoSum}(p_2, p_3)$
        \STATE $(s_1, t_1) := \mathrm{AVX2TwoSum}(q_1, p_4)$
        \STATE $s_2 := $\verb|_mm256_add_pd|($q_2, p_5$)
        \STATE $(s_1, t_0) := \mathrm{AVX2TwoSum}(s_1, t_0)$
        \STATE $s_2 :=$\verb|_mm256_add_pd|($s_2$, \verb|_mm256_add_pd|($t_0, t_1$))   
        \STATE $s_1 := $\verb|_mm256_add_pd|($s_1$, \verb|_mm256_mul_pd|(x[0], y[3]))
        \STATE $s_1 := $\verb|_mm256_add_pd|($s_1$, \verb|_mm256_mul_pd|(x[1], y[2]))
        \STATE $s_1 := $\verb|_mm256_add_pd|($s_1$, \verb|_mm256_mul_pd|(x[2], y[1]))
        \STATE $s_1 := $\verb|_mm256_add_pd|($s_1$, \verb|_mm256_mul_pd|(x[3], y[0]))
        \STATE $s_1 := $\verb|_mm256_add_pd|($s_1$, $q_0$)
        \STATE $s_1 := $\verb|_mm256_add_pd|($s_1$, $q_3$)
        \STATE $s_1 := $\verb|_mm256_add_pd|($s_1$, $q_4$)
        \STATE $s_1 := $\verb|_mm256_add_pd|($s_1$, $q_5$)
        \STATE $(r[0], r[1], r[2], r[3]) := \mathrm{AVX2Renorm}(p_0, p_1, s_0, s_1, s_2)$
        \STATE \textbf{return} ($r[0], r[1], r[2], r[3]$)
    \end{algorithmic}
\end{algorithm}

As later described, excluding AVX2Renorm, all parts of AVX2QDadd and AVX2QDmul were sufficient for a full implementation with AVX2; hence, a maximum speed up can be achieved with AVX2 rather than DD and TD. Additionally, they are completely faster than the MPFR 212-bits matrix multiplication. Therefore, the results obtained from the SIMDized QD arithmetic verify that the over quad-double MPF arithmetic with AVX2 is potentially more efficient than the MPFR arithmetic.

%
\subsection{DD, TD, and QD MFLOPS with element-wise addition and multiplication of vectors}

The vector data structure that we implemented for the multi-component MPF matrix multiplication is presented in this section. Although it is standard to use each MPF variable as one element of the array that expresses vector and matrix, we implemented one set of several arrays with each component of the vector or matrix element, as it is necessary for us to derive maximum performance using the \verb|_m256d| data-type of AVX2. Owing to these methods of expressing the MPF vector, highly performed continuous load and store instructions, such as \verb|__mm256_load_pd| and \verb|__mm256_store_pd| functions, can be utilized to read and store elements of vector or matrix by four elements at once.

\begin{figure}[htbp]
    \begin{center}
        \includegraphics[width=.6\textwidth]{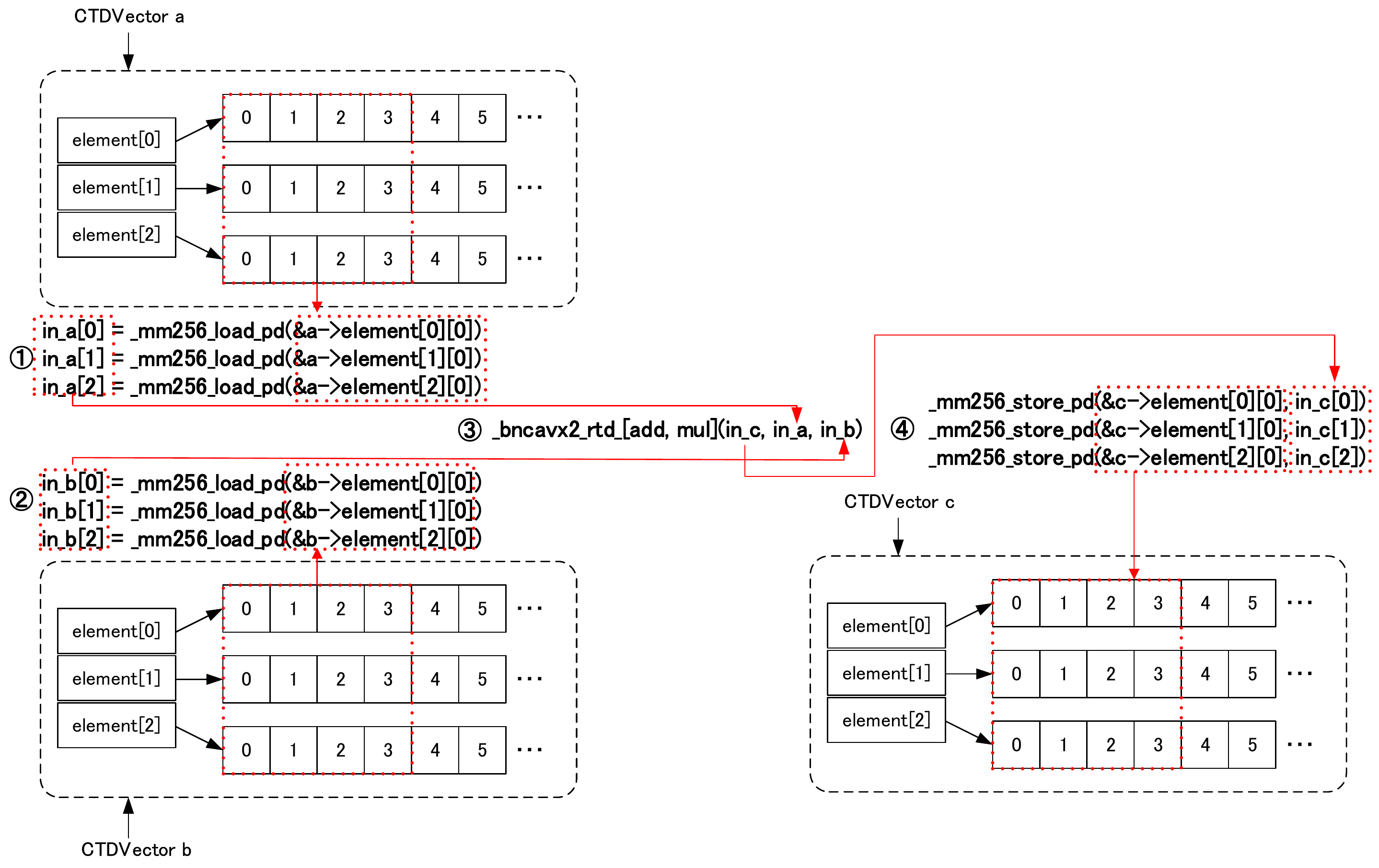}
        \caption{TD vector data type that can be treated with load/store instructions}\label{fig:simd_memory_access}
    \end{center}
\end{figure}

Therefore, our MPF vector and matrix data structure have two binary64, three binary64, and four binary64 arrays as DD, TD, and QD linear computations, respectively. By adopting these data structures, the TD linear computation case, as shown in \figurename\ref{fig:simd_memory_access}, two sets of three calls of load functions can settle down two operands of AVX2TD[add, mul] (\verb|bncavx2_rtd_[add, mul]|). Our adopted vector or matrix data structure for AVX2 can reduce the process time required to read and write as it eliminates the need to read each binary64 number.

In the remainder of this paper, we adopted two types of computational environments to evaluate our implemented codes.
\begin{description}
	\item[Corei9] Intel Core i9-10900X (3.6GHz, 10 cores), 16GB RAM Ubuntu 18.04.2, GCC 7.3.0
	\item[EPYC] AMD EPYC 7402P (1.5GHz, 24 cores), 64GB RAM, Ubuntu 18.04.5, GCC 7.5.0
\end{description}
%

\begin{figure}[htbp]
    \begin{center}
        \includegraphics[width=.35\textwidth]{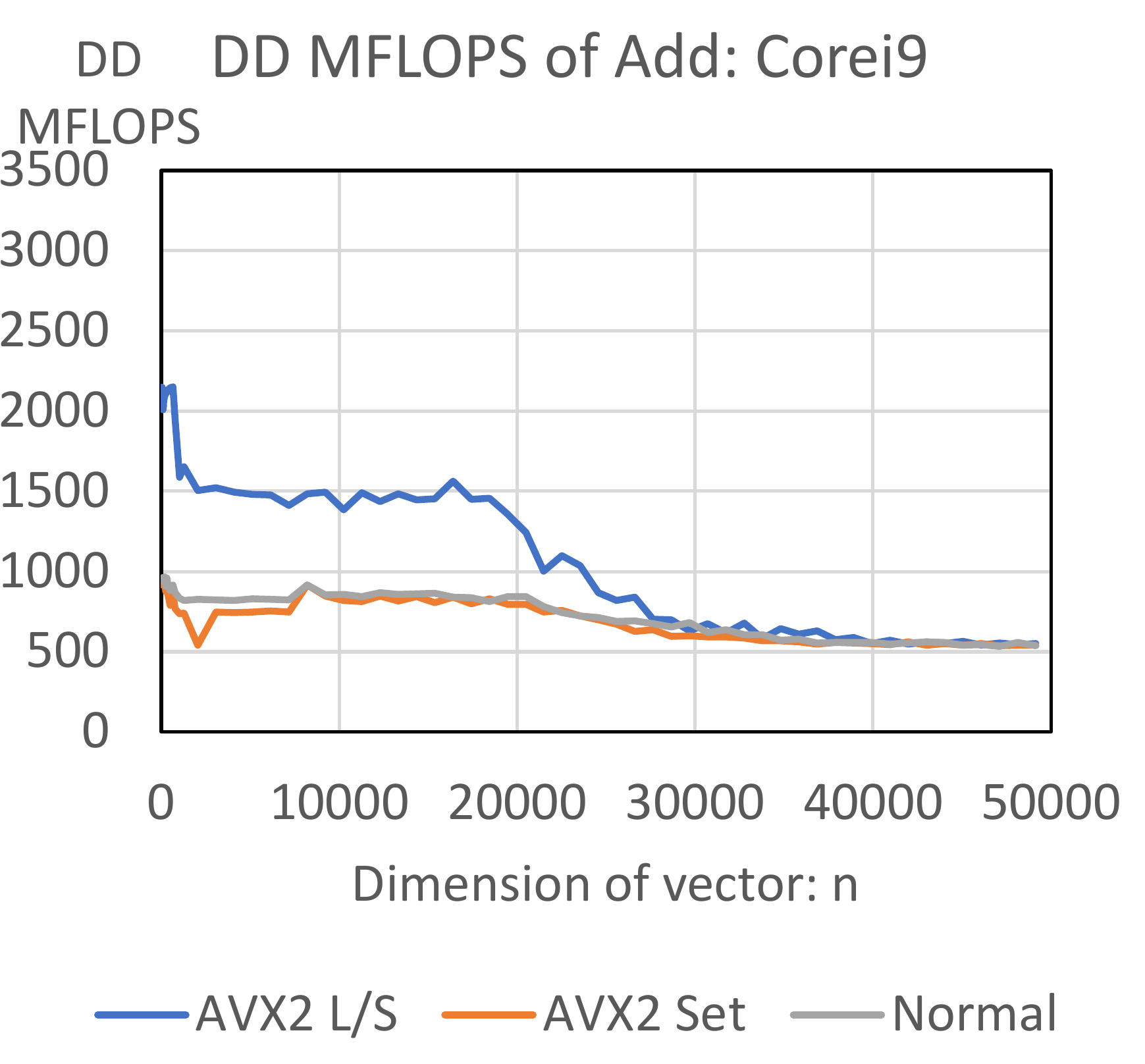}
        \includegraphics[width=.35\textwidth]{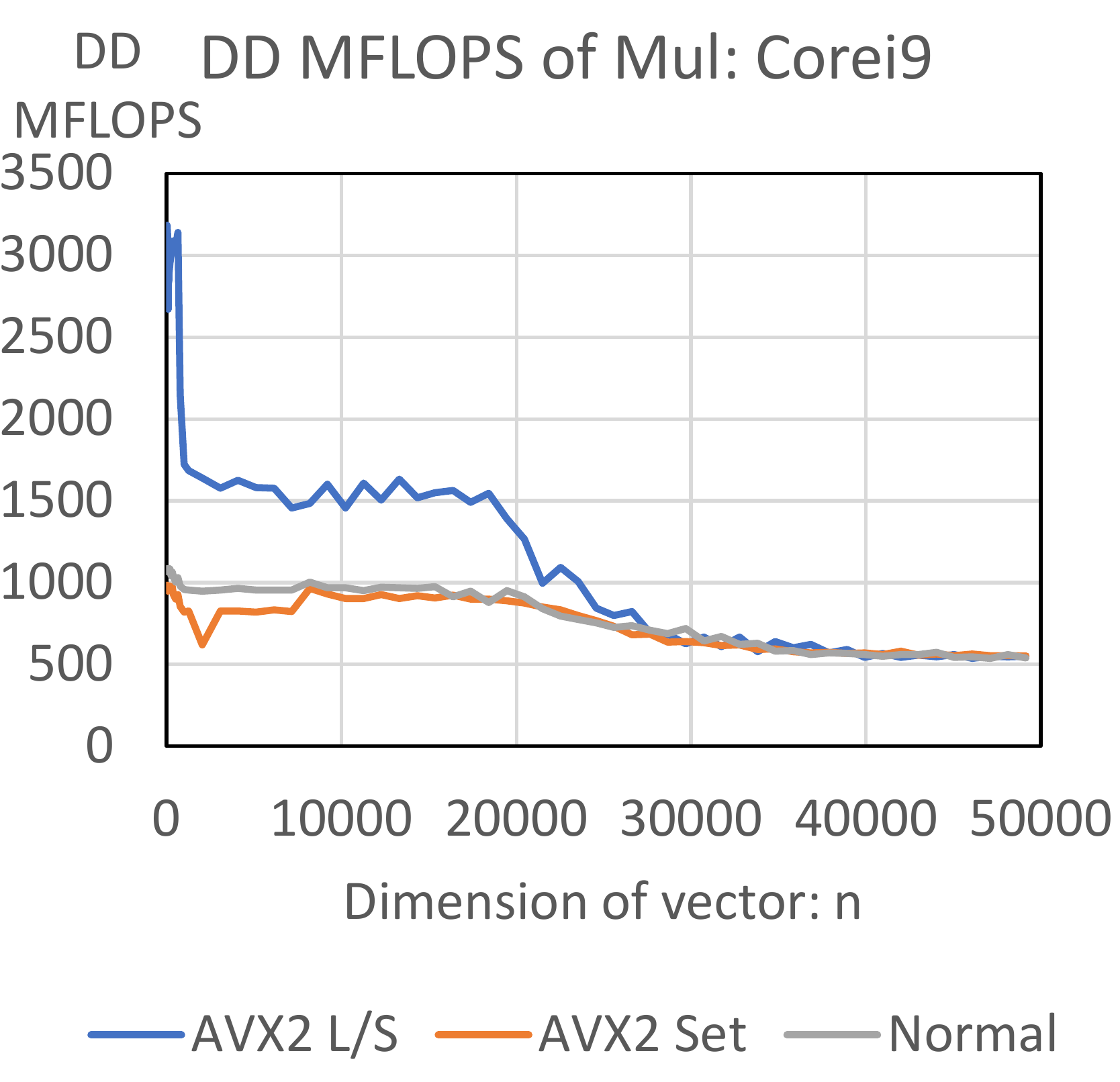}
        \caption{DD MFLOPS of DD element-wise addition and multiplication of vectors}\label{fig:ddmflops_addmul_corei9}
    \end{center}
\end{figure}

\begin{figure}[htbp]
    \begin{center}
        \includegraphics[width=.35\textwidth]{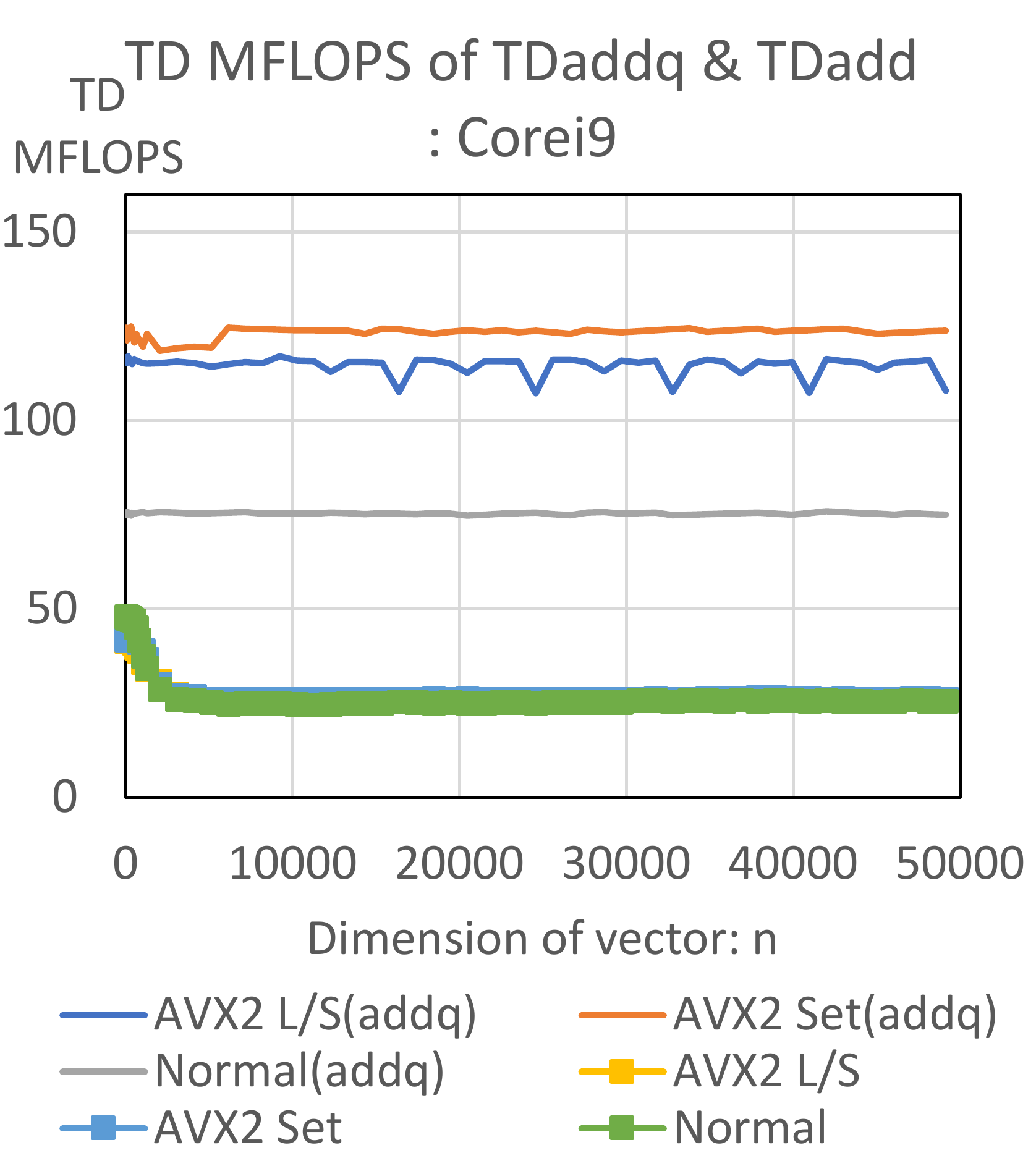}
        \includegraphics[width=.35\textwidth]{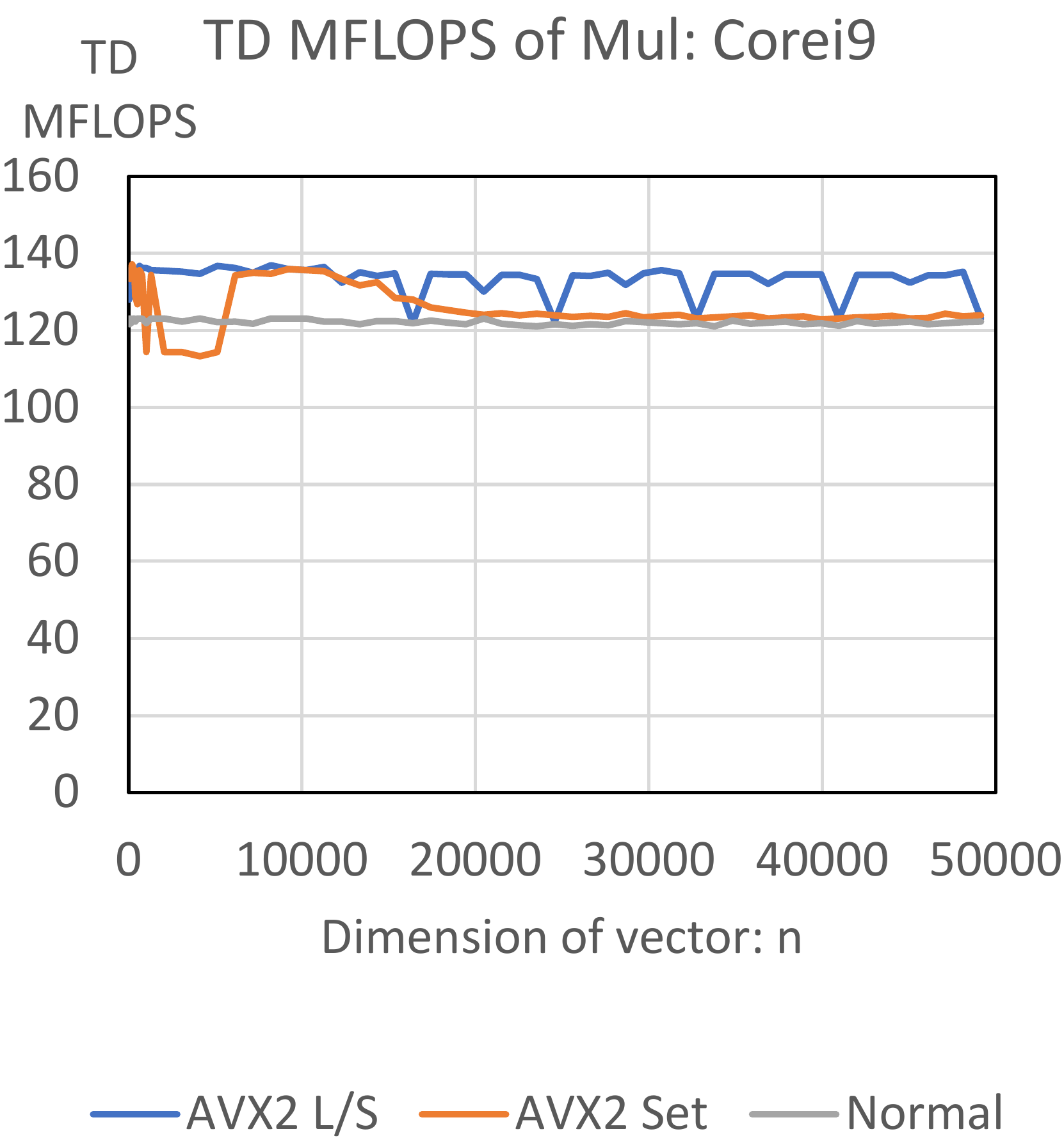}
        \caption{TD MFLOPS of TD element-wise addition and multiplication of vectors}\label{fig:tdmflops_addmul_corei9}
    \end{center}
\end{figure}

\begin{figure}[htbp]
    \begin{center}
        \includegraphics[width=.35\textwidth]{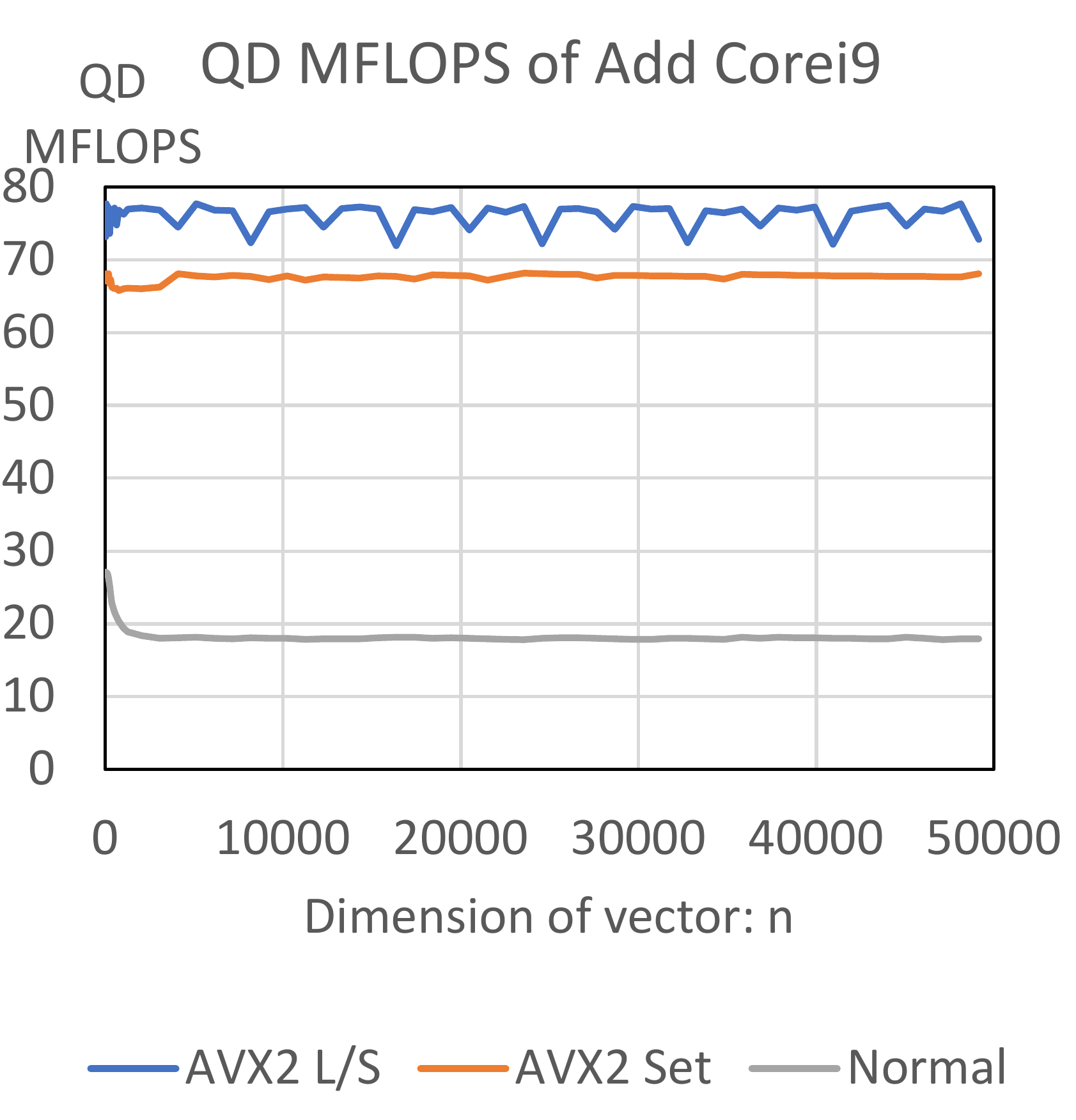}
        \includegraphics[width=.35\textwidth]{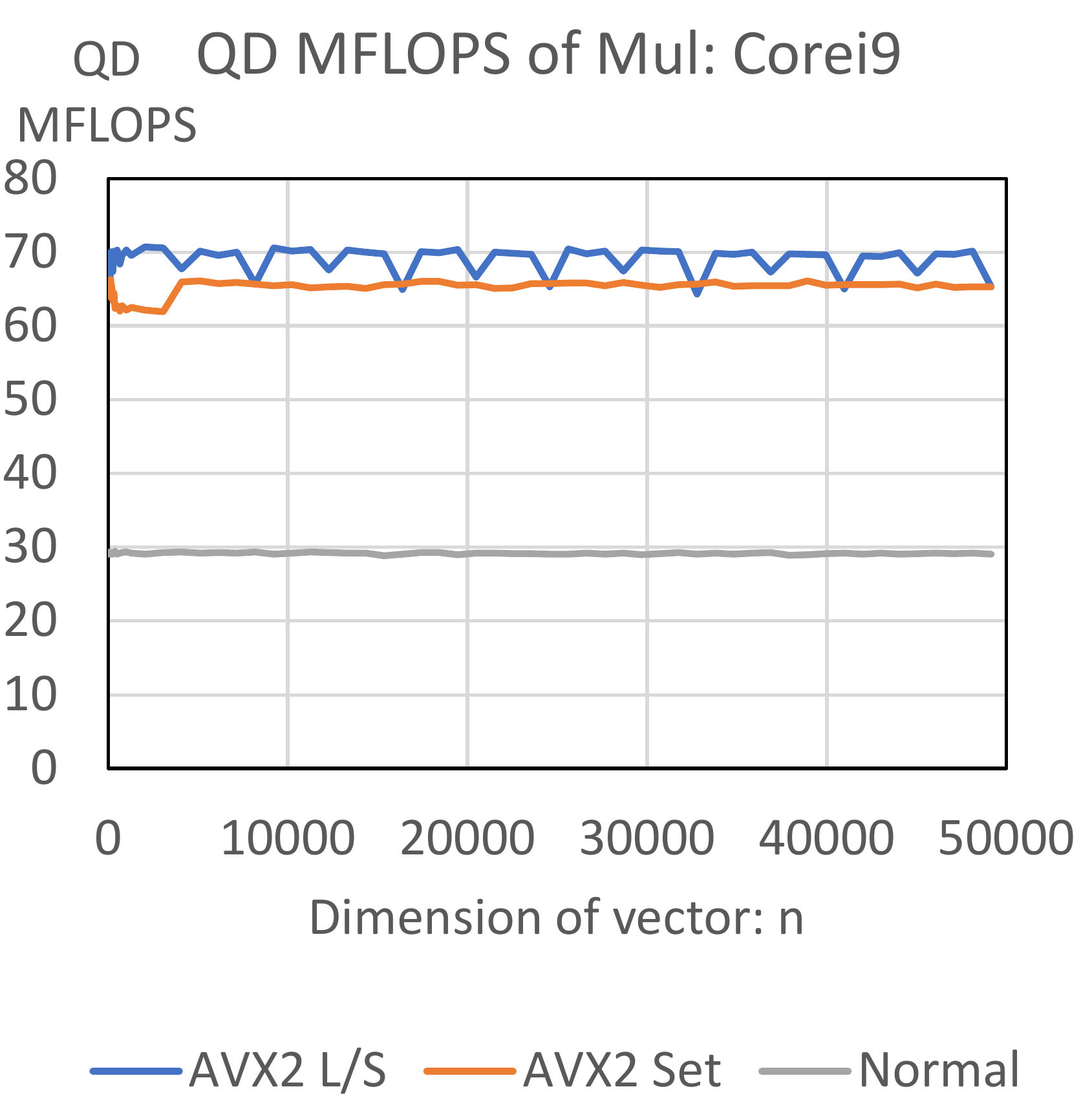}
        \caption{QD MFLOPS of QD element-wise addition and multiplication of vectors}\label{fig:qdmflops_addmul_corei9}
    \end{center}
\end{figure}

We demonstrated DD mega floating-point operations per second (MFLOPS), TD MFLOPS, and QD MFLOPS on a Corei9 environment to compare with three types of element-wise addition and multiplication of vectors in \figurename\ref{fig:ddmflops_addmul_corei9}, \figurename\ref{fig:tdmflops_addmul_corei9}, and \figurename\ref{fig:qdmflops_addmul_corei9}. In these FLOPS tests, we adopted two $n$-dimensional real vectors $\mathbf{a} = [a_1\ a_2\ ...\ a_n]^T$ and $\mathbf{b} = [b_1\ b_2\ ...\ b_n]^T$ as random numbers. Thereafter, we evaluated the number of additions or multiplications per second by computing each element, such as $a_i + b_i$ or $a_i\cdot b_i$. ``AVX2 L/S'' depicts MFLOPS with SIMDized calculation and AVX2 Load/Store instructions, ``AVX2 Set'' SIMDized calculation and AVX2 Set and normal substitution to arrays, and complete ``Normal'' standard calculation, except for AVX2 instructions. The DD, TD, and QD MFLOPS values are presented in the vertical axis.

\figurename\ref{fig:ddmflops_addmul_corei9} illustrates three types of DD MFLOPS. DD arithmetic requires less computation than TD or QD; therefore, the utilization of AVX2 load/store instructions is largely effective for the cache memory of a CPU. In fact, when all elements of the vectors are stored in cache memory, the DD MFLOPS is 1.5 times larger than the other two types of DD computations. For vector sizes beyond cache memory, both DD addition and multiplication for all types of calculations have the same performance.

\figurename\ref{fig:tdmflops_addmul_corei9} presents the TD MFLOPS for different processes. As can be observed in the MFLOPS of the original triple-double addition (TDadd and AVX2TDadd) presented in the left figure, they realized approximately 26 TD MFLOPS owing to the poor performance of the merge function in Algorithm \ref{algo:AVX2TDsum}; hence, SIMDization is ineffective for them. In contrast, TDaddq performed well at 75 TD MFLOPS, and AVX2TDaddq with load/store instructions also achieved approximately 115 TD MFLOPS, with AVX2 set at 123 TD MFLOPS. The Corei9 environment is the only case in which the AVX2 load/store performance is less than that of the AVX2 set. Otherwise, TDmul cannot achieve sufficient efficiency with AVX2, as more than approximately 20 TD MFLOPS is achieved. Through these performance tests, we can infer that cache memory is not effective for TD arithmetic.

\figurename\ref{fig:qdmflops_addmul_corei9} illustrates QD MFLOPS on the Corei9 environment. Similar to the TD arithmetic, cache memory is never effective for QD arithmetic, and AVX2 can provide maximum performance for both QD addition and multiplication. AVX2QDadd is four times faster than normal QDadd, whereas AVX2QDmul is two times faster than normal QDmul. We can expect better performance for longer multi-component MPF arithmetic with AVX2.

As demonstrated above, these MPF MFLOPS and SIMDized arithmetic with AVX2 load/store instructions can achieve the best performance for almost all cases. These results are similar to those of EPYC. Therefore, we can expect high-performance matrix multiplication with AVX2 for DD, TD, and QD arithmetic.

%
\section{Benchmark tests of DD, TD, and QD matrix multiplication}

To evaluate the performance of our multi-component MPF matrix multiplication, we considered a dense real matrix multiplication to obtain $C := AB\in\mathbb{R}^{n\times n}$, where $A =[a_{ij}]$ $\in \mathbb{R}^{n\times n}$, $B = [b_{ij}]$$\in \mathbb{R}^{n\times n}$ and $c_{ij}$, the elements of $C$ are obtained as 
\begin{equation}
	c_{ij} := \sum^n_{k = 1} a_{ik} b_{kj}, \label{eqn:matrix_mul_simple}
\end{equation}

To draw the power of cache memories in the majority of current CPUs, ``block'' algorithms are adopted in dense matrix multiplications, which normally divide $A$ and $B$ into small block matrices with less rows and columns than the block size $n_{ \rm min}$. We fixed $n_{\rm min} = 32$ for our benchmark tests in this paper.

\begin{figure}[htbp]
    \begin{center}
        \includegraphics[width=.4\textwidth]{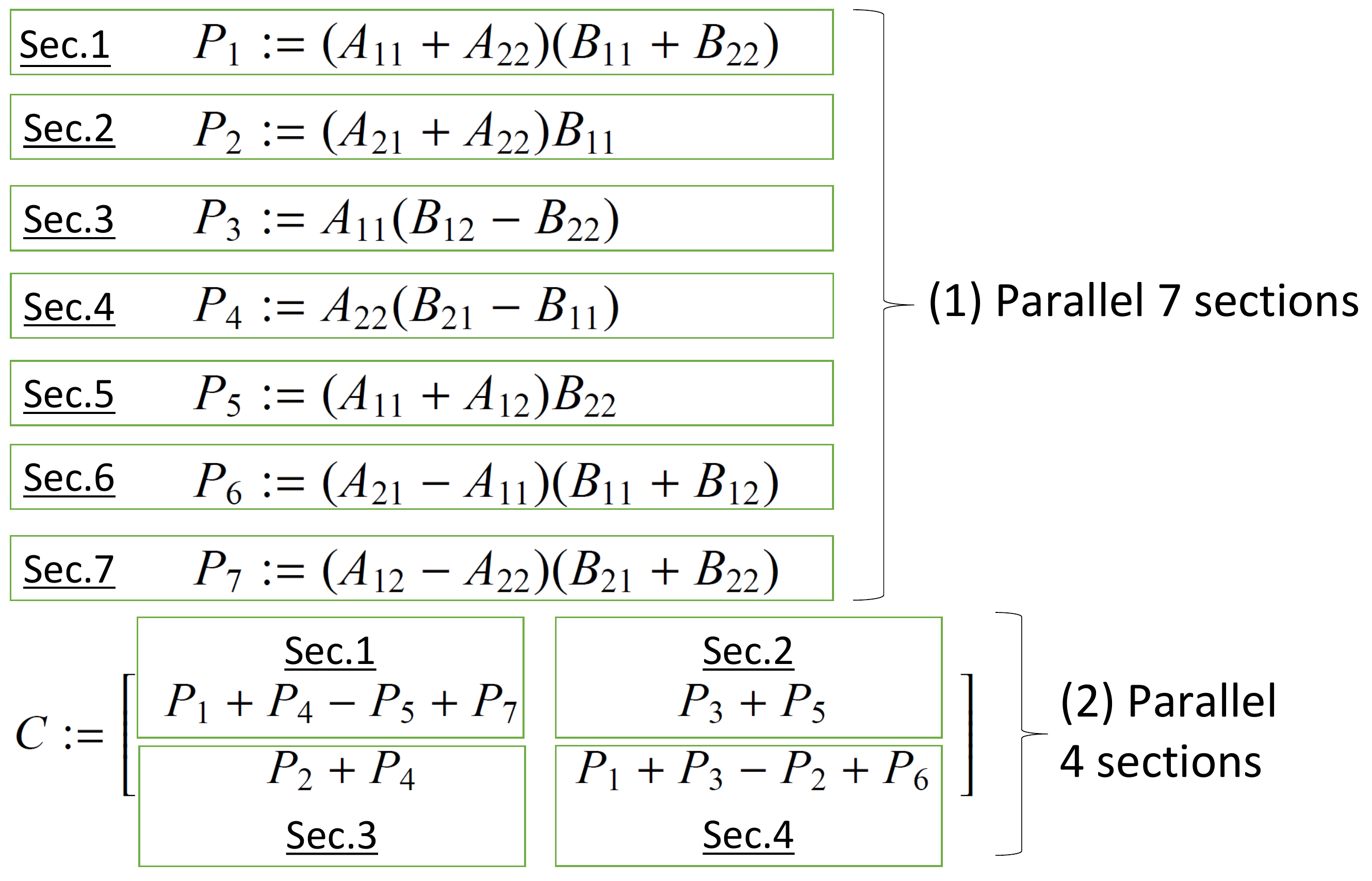}
        \caption{Strassen matrix multiplication and its parallelization}\label{fig:matrix_mul_algorithms}
    \end{center}
\end{figure}

The ``Strassen'' algorithm is a  divide-and-conquer method. For even dimensional matrices $A$ and $B$, four divided matrices are prepared and $C := AB$ is obtained, as shown in \figurename \ref{fig:matrix_mul_algorithms}.

Among all types of matrix multiplication, we adopted the row-major method to access the elements of matrices and parallelized them by OpenMP\cite{kouya_strassen2016}.

%

We fixed matrices $A$ and $B$ as follows:
\[ A = \left[\sqrt{5} \left(i + j - 1\right)\right]^{n}_{i,j=1},\ B = \left[\sqrt{3}\left(n - i\right)\right]^{n}_{i,j=1}. \]
Since all elements were positive, there was no significant loss of significant digits in the matrix multiplication. We determined one or two decimal digits of loss for all types of matrix multiplications.

All the codes for matrix multiplication were written in C or C++ languages, and the major programs adopted for the benchmark, including Rgemm of MPBLAS, were compiled with options as follows:
\begin{verbatim}
g++ -O3 -std=c++11 -mavx2 -mfma -fopenmp
\end{verbatim}

MPBLAS, which was downloaded on June 2019 from Github, was installed and used in our C++ benchmark programs.

%
\subsection{Computational time of serial matrix multiplication and comparison with Rgemm (MPBLAS)}

First, we present the results obtained by serial matrix multiplication on both the Corei9 and EPYC environments in \tablename\ \ref{table:serial_table}. In the tables, the ``B'' fields, ``B + A,'' ``S'' fields, and ``S + A'' represent the computational time (s) of block matrix multiplication, block matrix multiplication with AVX2, Strassen matrix multiplication, and Strassen matrix multiplication with AVX2, respectively. The computational time (``M'' fields) provided by the Rgemm of MPBLAS are presented in the rightmost columns of the tables. The underlined values in the tables below mean minimal computational time for the number of dimension $n$ on the environment.

We have previously reported that the serial Strassen matrix multiplication and parallelized versions for the large $n$ are more efficient than Rgemm of MPBLAS \cite{kouya_strassen2016}. Our accelerated versions of block and Strassen matrix multiplication with AVX2 can be achieved more than two times faster than the MPBLAS for any size of matrices. Although the TD arithmetic is currently not supported by MPBLAS, we can confirm that the computational time of TD matrix multiplication is located between DD and QD arithmetic for any size and any type of algorithm.

{\tiny \begin{table}[htbp]
\begin{center}\small
\caption{Serial computational time of matrix multiplication (unit: s): Corei9(left) and EPYC(right)}\label{table:serial_table}
{\small \begin{tabular}{c|cc|cc|c}
\multicolumn{6}{c}{DD : Corei9} \\ \hline
$n$  & B  & B+A & S & S+A & M \\ \hline
1023 & 7.84 & 2.46 & 4.35 & \underline{1.57} & 5.98 \\
1024 & 7.86 & 2.46 & 4.34 & \underline{1.55} & 6.00 \\
1025 & 8.61 & 2.68 & 4.40 & \underline{1.59} & 6.01 \\ \hline
4095 & 507.8 & 162.7 & 212.3 & \underline{74.47} & 390.7 \\
4096 & 509.1 & 161.8 & 212.5 & \underline{74.07} & 390.7 \\
4097 & 518.0 & 161.7 & 213.3 & \underline{74.94} & 391.0 \\ \hline
\multicolumn{6}{c}{TD : Corei9} \\ \hline
1023 &	50.06	&20.54	&26.75	&\underline{12.64} & N/A \\
1024 &	50.10	&20.54	&26.96	&\underline{12.48} & N/A \\
1025 &	54.28	&22.56	&26.95	&\underline{12.61} & N/A \\ \hline
4095 &	3202	&1316	&1276	&\underline{619.5} & N/A \\
4096 &	3205	&1317	&1276	&\underline{618.6} & N/A \\
4097 &	3272	&1345	&1287	&\underline{620.5} & N/A \\ \hline
\multicolumn{6}{c}{QD : Corei9} \\ \hline
1023 & 102.8 & 31.41 & 54.90 & \underline{19.55} & 73.76 \\
1024 & 102.9 & 31.40 & 54.71 & \underline{19.41} & 76.04 \\
1025 & 111.8 & 34.41 & 55.23 & \underline{19.63} & 74.12 \\ \hline
4095 & 6491 & 2013 & 2720 & \underline{970.5} & 4729 \\
4096 & 6493 & 2013 & 2721 & \underline{968.6} & 4730 \\
4097 & 6624 & 2059 & 2725 & \underline{972.6} & 4727 \\ \hline
\end{tabular}}
{\small \begin{tabular}{c|cc|cc|c}
\multicolumn{6}{c}{DD : EPYC} \\ \hline
$n$ & B  & B+A & S & S+A & M \\ \hline
1023 & 9.27  & 2.37  & 5.00  & \underline{1.43}  & 7.55 \\
1024 & 9.31  & 2.37  & 4.98  & \underline{1.41}  & 7.54 \\
1025 & 10.16  & 2.57  & 5.03  & \underline{1.44}  & 7.59 \\ \hline
4095 & 595.3  & 163.6  & 243.6  & \underline{68.04}  & 482.8 \\
4096 & 609.1  & 163.6  & 243.6  & \underline{67.75}  & 484.0 \\
4097 & 611.4  & 152.6  & 244.7  & \underline{71.03}  & 483.8 \\ \hline
\multicolumn{6}{c}{TD : EPYC} \\ \hline
1023 &	61.12	&21.19	&32.18	&\underline{12.97} & N/A \\ 
1024 &	61.15	&21.16	&31.98	&\underline{12.78} & N/A \\
1025 &	66.46	&23.17	&32.20	&\underline{12.89} & N/A \\ \hline
4095 &	3918	&1379	&1555	&\underline{632.3} & N/A \\ 
4096 &	3934	&1378	&1554	&\underline{631.5} & N/A \\
4097 &	4002	&1387	&1561	&\underline{634.2} & N/A \\ \hline
\multicolumn{6}{c}{QD : EPYC} \\ \hline
1023 & 127.8 & 42.82 & 71.23 & \underline{25.17} & 83.95 \\
1024 & 127.8 & 42.78 & 71.10 & \underline{25.07} & 84.18 \\
1025 & 139.6 & 46.88 & 71.55 & \underline{25.32} & 84.47 \\ \hline
4095 & 8177 & 2753 & 3440 & \underline{1245} & 5384 \\
4096 & 8188 & 2754 & 3438 & \underline{1243} & 5390 \\
4097 & 8371 & 2808 & 3449 & \underline{1247} & 5394 \\ \hline
\end{tabular}}
\end{center}
\end{table}}
%
\subsection{Speedup ratio of parallelization and computational time}

\begin{figure}[htbp]
    \begin{center}
        \includegraphics[width=.35\textwidth]{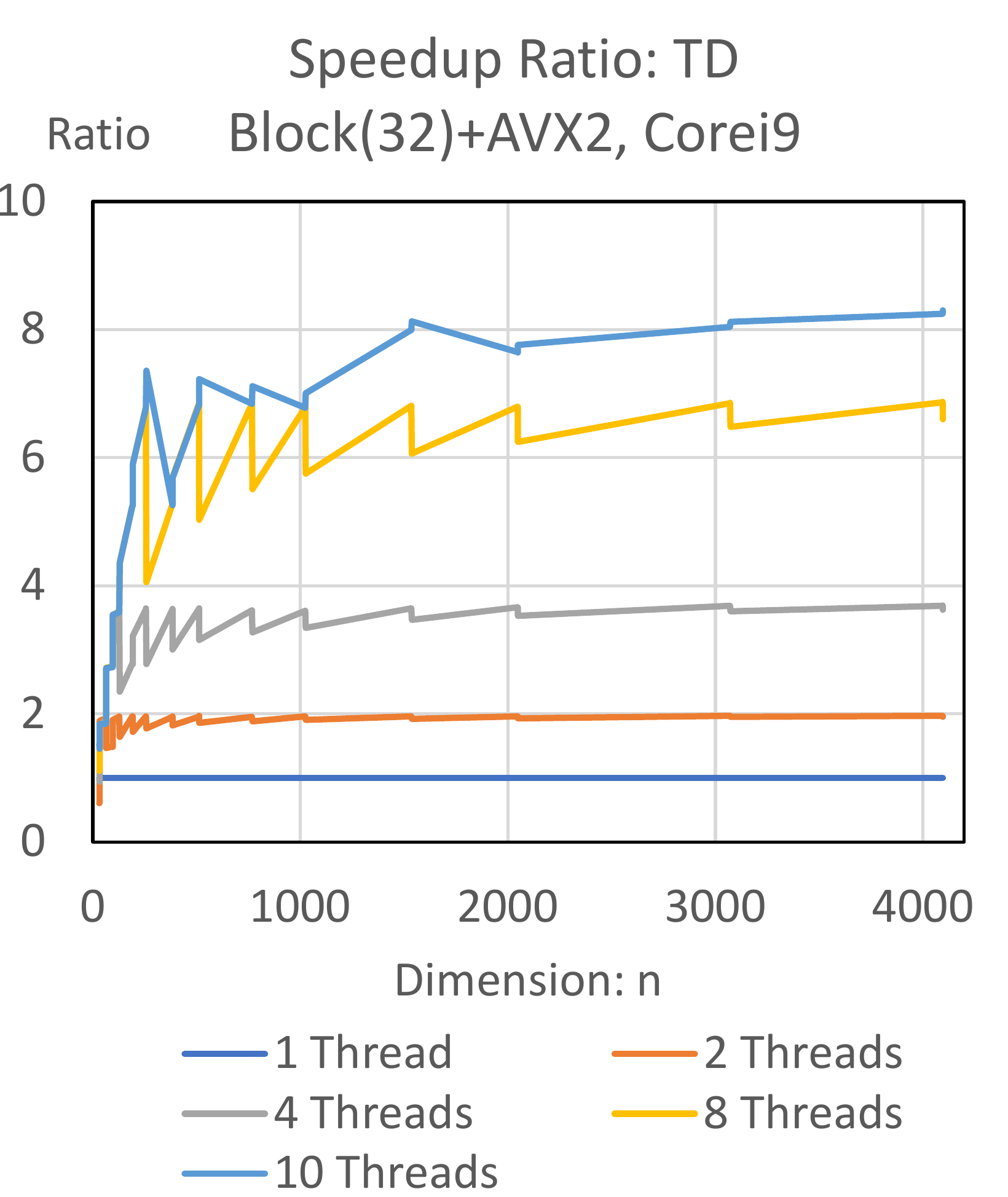}
        \includegraphics[width=.35\textwidth]{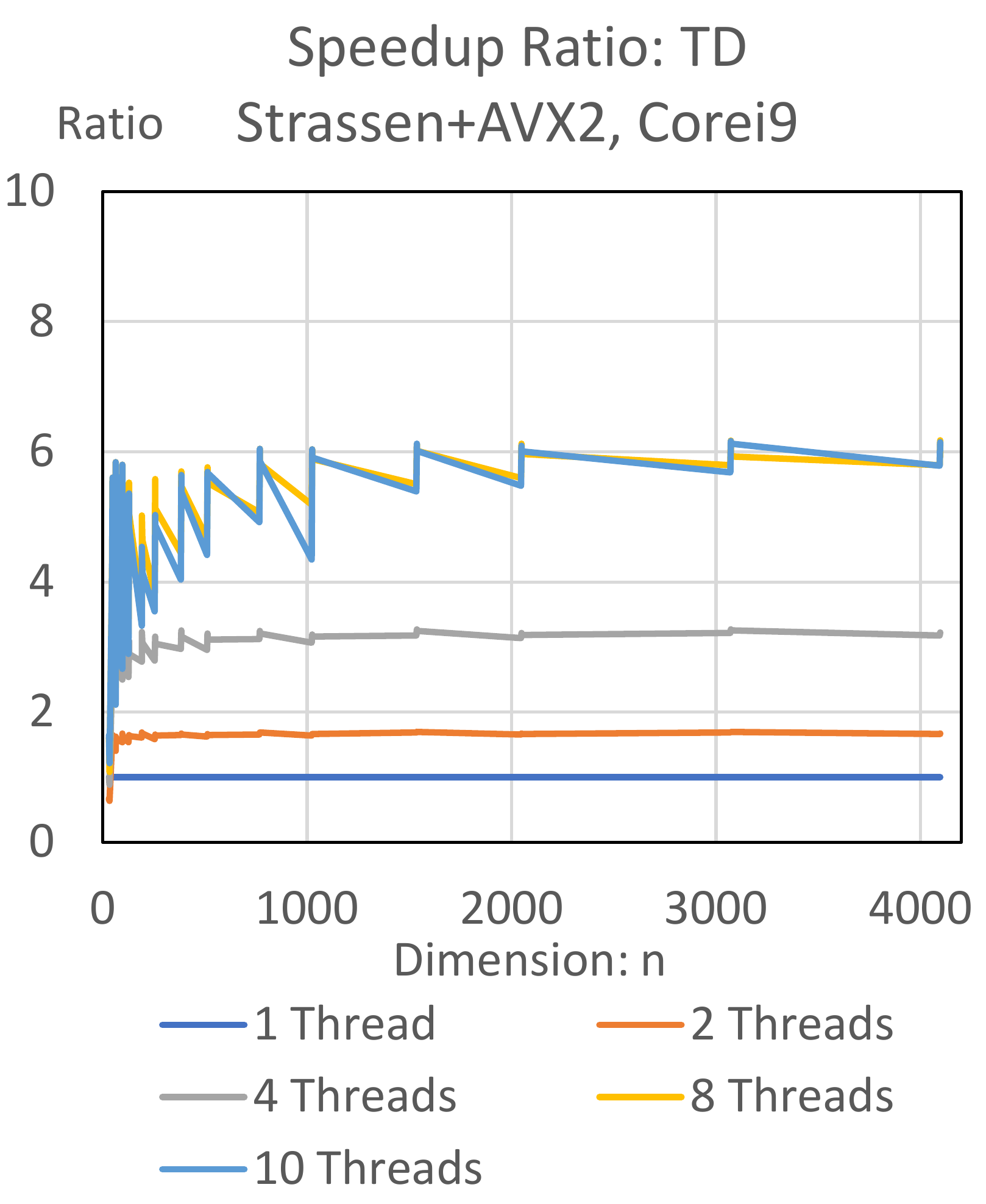}
        \caption{Speedup ratio of paralyzed TD block matrix multiplication on Corei9}\label{fig:speedup_td_corei9}
    \end{center}
\end{figure}

Next, we present the results of parallelization using OpenMP. As later explained, it is unusual to perform parallel execution with over 16 threads on EPYC. As a typical case of the speedup ratio of parallelization, we presented \figurename\ref{fig:speedup_td_corei9} on Corei9. As shown in this figure, we can infer that stable speedup has been achieved for block matrix multiplication by approximately 10 threads and that the speedup ratio of the Strassen matrix multiplication has been limited by approximately six times using more than eight threads. This limitation occurs owing to Strassen’s method of parallelization, which produces new threads in recursive calls. We assume that the scheduling of many threads in the calls is not smoothly executed. 

\begin{figure}[htbp]
    \begin{center}
        \includegraphics[width=.35\textwidth]{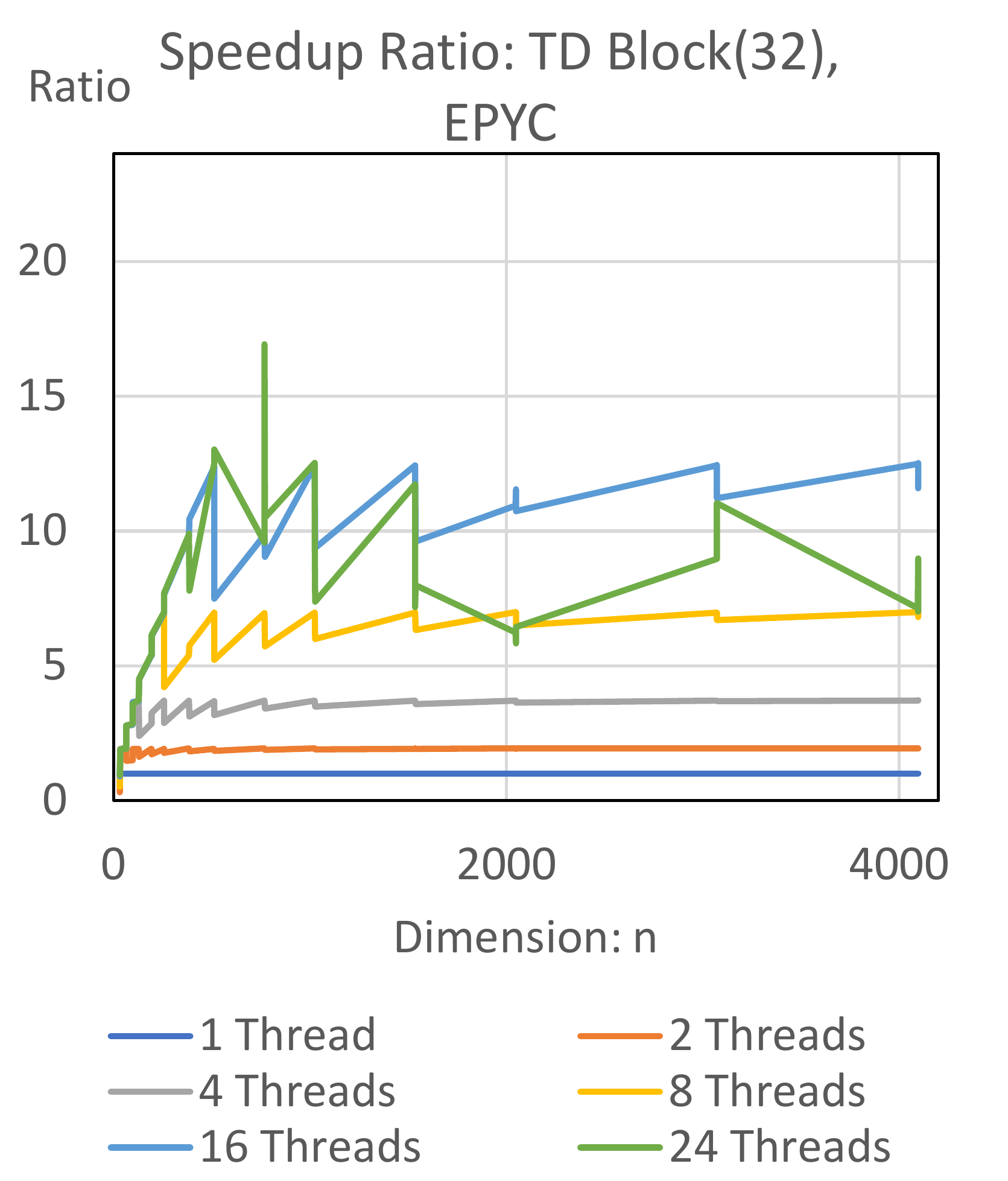}
        \includegraphics[width=.35\textwidth]{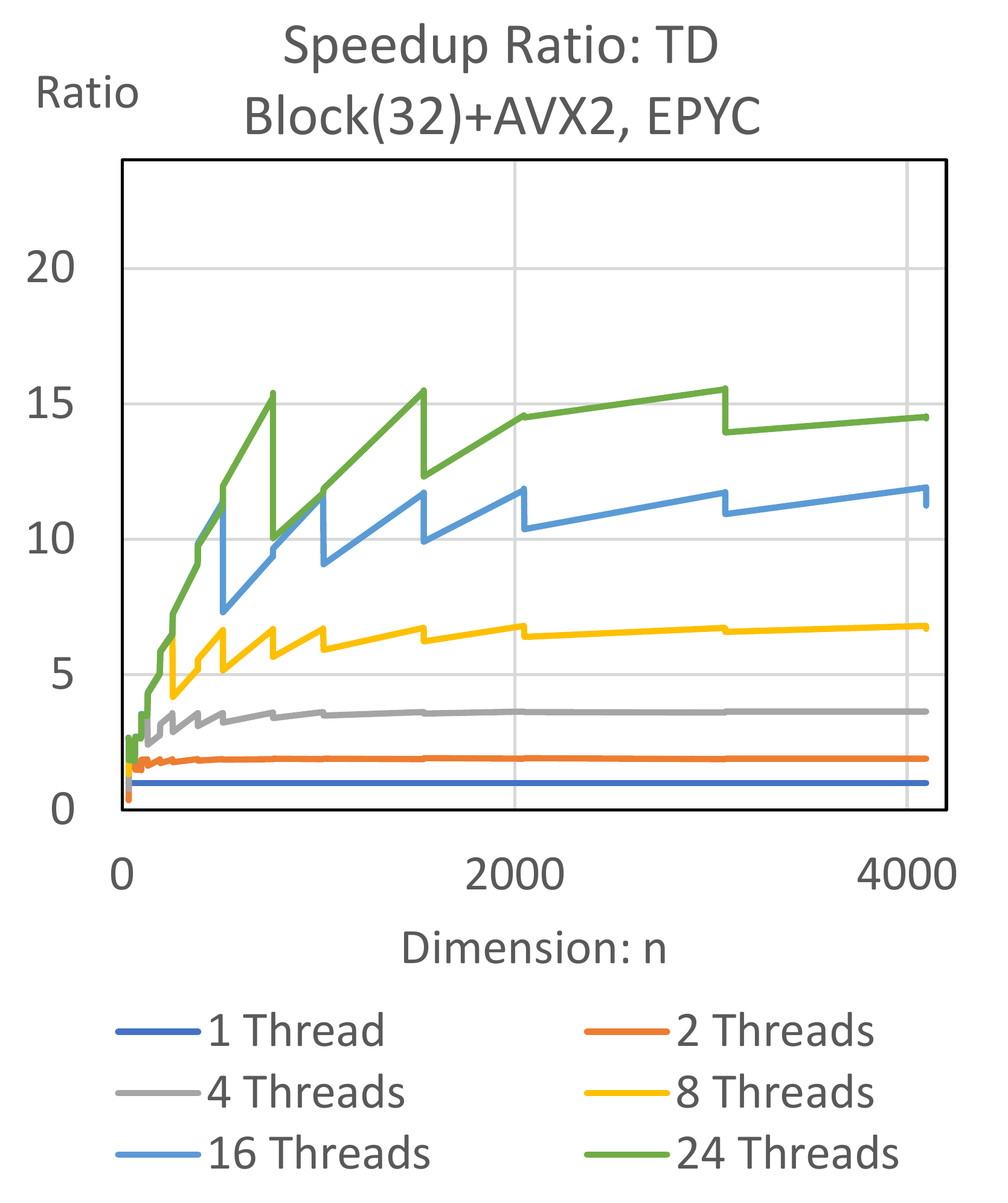}
        \caption{Speedup ratio of parallelized TD block matrix multiplication on EPYC}\label{fig:speedup_td_epyc}
    \end{center}
\end{figure}

In contrast to the Corei9 environment, although the EPYC environment always suspends the performance of parallelization with over 16 threads for DD, TD, and QD matrix multiplications (the left figure of \figurename\ref{fig:speedup_td_epyc}), AVX2 can eliminate the limitation of parallelization, as demonstrated in the right side of \figurename\ref{fig:speedup_td_epyc}. The DD and QD matrix multiplication on EPYC can be described as TD.

The limitation of the speedup ratio of the parallelized Strassen matrix multiplication is demonstrated similar to that in the Corei9 environment, and AVX2 cannot eliminate the failure of parallelization. We infer that to improve the performance of the parallelized block matrix multiplication for over eight threads on CPUs with many cores, it is necessary to completely rewrite the Strassen routine.

\tablename\ \ref{table:parallel_table} presents all computational time (s) with the same number of threads similar to those of the cores on Corei9 (10 cores) and EPYC (24 cores).

{\tiny \begin{table}[htbp]
\begin{center}\small
\caption{Seconds of parallelized matrix multiplication: 10 threads on Corei9 (left) and 24 threads on EPYC (right)}\label{table:parallel_table}
{\small \begin{tabular}{c|cc|cc|c}
\multicolumn{5}{c}{DD : Corei9 10 Threads} \\ \hline
$n$  & B  & B+A & S & S+A \\ \hline
1023 & 1.15 & \underline{0.35} & 0.86 & 0.39 \\
1024 & 1.15 & 0.35 & 0.69 & \underline{0.33} \\
1025 & 1.23 & \underline{0.37} & 0.79 & 0.38 \\ \hline
4095 & 61.51 & 19.68 & 33.28 & \underline{16.63} \\
4096 & 61.97 & 19.67 & 32.60 & \underline{16.19} \\
4097 & 62.65 & 19.42 & 34.20 & \underline{16.44} \\ \hline
\multicolumn{5}{c}{TD : Corei9 10 Threads} \\ \hline
1023 & 7.65 & 3.03 & 4.57 & \underline{2.91} \\
1024 & 7.67 & 3.02 & 4.11 & \underline{2.07} \\
1025 & 8.12 & 3.22 & 4.30 & \underline{2.13} \\ \hline
4095 & 409.94 & 159.68 & 202.34 & \underline{107.08} \\
4096 & 410.48 & 159.66 & 196.19 & \underline{100.55} \\
4097 & 416.29 & 162.04 & 198.99 & \underline{101.92} \\ \hline
\multicolumn{5}{c}{QD : Corei9 10 Threads} \\ \hline
1023 & 15.23 & 4.65 & 9.12 & \underline{3.43} \\
1024 & 15.22 & 4.64 & 8.73 & \underline{3.13} \\
1025 & 16.18 & 4.95 & 8.88 & \underline{3.20} \\ \hline
4095 & 824.23 & 247.43 & 437.68 & \underline{158.57} \\
4096 & 822.68 & 247.48 & 431.15 & \underline{153.63} \\
4097 & 834.92 & 251.47 & 434.29 & \underline{155.50} \\ \hline
\end{tabular}} {\small \begin{tabular}{c|cc|cc}
\multicolumn{5}{c}{DD : EPYC 24 Threads} \\ \hline
$n$  & B  & B+A & S & S+A \\ \hline
1023 & 0.78 & \underline{0.20} & 0.80 & 0.59 \\
1024 & 0.73 & \underline{0.20} & 0.77 & 0.54 \\
1025 & 0.79 & \underline{0.22} & 0.84 & 0.57 \\ \hline
4095 & 160.68 & \underline{10.74} & 35.67 & 25.43 \\
4096 & 191.24 & \underline{10.72} & 35.24 & 24.90 \\
4097 & 169.09 & \underline{10.28} & 36.26 & 25.23 \\ \hline
\multicolumn{5}{c}{TD : EPYC 24 Threads} \\ \hline
1023 & 4.88 & \underline{1.81} & 5.16 & 2.80 \\
1024 & 6.18 & \underline{1.81} & 4.68 & 2.30 \\
1025 & 9.01 & \underline{1.95} & 4.85 & 2.36 \\ \hline
4095 & 549.21 & \underline{95.02} & 229.31 & 117.44 \\
4096 & 437.79 & \underline{94.76} & 223.25 & 109.44 \\
4097 & 570.95 & \underline{95.98} & 225.54 & 110.36 \\ \hline
\multicolumn{5}{c}{QD : EPYC 24 Threads} \\ \hline
1023 & 26.95 & \underline{3.35} & 8.91 & 4.07 \\
1024 & 28.98 & \underline{3.37} & 8.40 & 3.50 \\
1025 & 27.35 & 3.60 & 8.53 & \underline{3.56} \\ \hline
4095 & 863.46 & \underline{171.27} & 415.35 & 179.00 \\
4096 & 865.12 & 171.13 & 407.80 & \underline{170.27} \\
4097 & 882.51 & 172.95 & 410.26 & \underline{171.72} \\ \hline
\end{tabular}}
\end{center}
\end{table}}

As shown in \tablename\ \ref{table:parallel_table}, block matrix multiplication tends to be increasingly faster than the Strassen on a CPU with many cores, owing to the limitations posed by parallelization of Strassen matrix multiplication. In fact, Strassen matrix multiplication is the fastest on Corei9, and more cases of shorter time of block matrix multiplication are observed on EPYC. We expect parallelized block matrix multiplication with DD, TD, and QD precision to be the fastest on a CPU with over 32 cores.

%
\section{Conclusion and future work}

Based on the accelerated DD, TD, and QD MPF matrix multiplication implementations in AVX2, we can confirm that block matrix and Strassen matrix multiplication show much higher speedup than the previous implementations. Furthermore, we can confirm that the speedup ratio is not limited by paralleling it with OpenMP, and it is stabilized in the EPYC environment. 

Although the current implementation of parallelized Strassen matrix multiplication has exhibited a limitation of acceleration in CPUs with more than eight cores, BLAS level 1 and 2 functions have been implemented owing to their required implementations in Strassen matrix multiplication. We will apply and develop linear and nonlinear algorithms with our accelerated MPF linear computation library in our future work.

%
\section*{Acknowledgment}
This work was supported by JSPS KAKENHI Grant Number JP20K11843.

%

\end{document}